\documentstyle{amlts}
\begin{document}
\annalsline{155}{2002}
\received{January 25, 2000}
\startingpage{105}
\def\bye{\end{document}}
 \font\tenrm=cmr10
\input boxedeps.tex 
\SetepsfEPSFSpecial 
\HideDisplacementBoxes
\def\figin#1#2{ 
$$
 {\BoxedEPSF{#1.eps scaled
#2}%
}%
$$
\noindent}
\catcode`\@=11
\font\twelvemsb=msbm10 scaled 1100
\font\tenmsb=msbm10
\font\ninemsb=msbm10 scaled 800
\newfam\msbfam
\textfont\msbfam=\twelvemsb  \scriptfont\msbfam=\ninemsb
  \scriptscriptfont\msbfam=\ninemsb
\def\msb@{\hexnumber@\msbfam}
\def\Bbb{\relax\ifmmode\let\next\Bbb@\else
 \def\next{\errmessage{Use \string\Bbb\space only in math
mode}}\fi\next}
\def\Bbb@#1{{\Bbb@@{#1}}}
\def\Bbb@@#1{\fam\msbfam#1}
\catcode`\@=12

 \catcode`\@=11
\font\twelveeuf=eufm10 scaled 1100
\font\teneuf=eufm10
\font\nineeuf=eufm7 scaled 1100
\newfam\euffam
\textfont\euffam=\twelveeuf  \scriptfont\euffam=\teneuf
  \scriptscriptfont\euffam=\nineeuf
\def\euf@{\hexnumber@\euffam}
\def\frak{\relax\ifmmode\let\next\frak@\else
 \def\next{\errmessage{Use \string\frak\space only in math
mode}}\fi\next}
\def\frak@#1{{\frak@@{#1}}}
\def\frak@@#1{\fam\euffam#1}
\catcode`\@=12

\newcommand{\binom}[2]{\left(\begin{array}{c}#1\\#2\end{array}\right)}
\def\C{{\Bbb C}}
\def\Z{{\Bbb Z}}
\def\bC{{\overline{\Bbb C}}}
\def\R{{\Bbb R}}
\def\T{{\Bbb T}}
\def\U{{\Bbb U}}
\def\bU{{ \overline{\Bbb U}}}
\def\card{{{\rm card\,}}}
\def\id{{{\rm id}}}
\def\bd{{{\rm bd\,}}}
\def\diam{{{\rm diam\,}}}
\def\dim{{{\rm dim\,}}}
\def\dist{{{\rm dist\,}}}
\def\Poly{{{\rm Poly}}}
\def\s{{{\bf s}}}
\def\CP{{\Bbb CP}}
\def\RP{{\Bbb RP}}
\def\SS{{\Bbb S}}
\def\Tilde{\widetilde}

 \title{Rational functions with real critical points\\ and the B. and M. Shapiro
conjecture\\ in real enumerative geometry} 
\shorttitle{Rational functions with real critical points} 

  \acknowledgements{The first author was supported by NSF grant DMS-0100512
and by Bar Ilan University.  The second author was supported by NSF grant DMS-0070666 and by MSRI.}
 \twoauthors{A. Eremenko}{A. Gabrielov}
 \institutions{Purdue University, West Lafayette, IN\\
{\eightpoint {\it E-mail addresses\/}: eremenko@math.purdue.edu}\\
\hglue.97in {\eightpoint agabriel@math.purdue.edu}} 


 \bigbreak
\centerline{\bf Abstract}
\vglue12pt

Suppose that $2d-2$ tangent lines to the rational normal curve
$z\mapsto\break (1:z:\ldots:z^d)$ in $d$-dimensional complex projective space are given.
It was known that the number of codimension $2$ subspaces intersecting all these
lines is always finite; for a generic configuration it is equal to the $d^{\rm th}$ Catalan number.
We prove that for real tangent lines, all these codimension $2$ subspaces are
also real, thus confirming a special case of a general conjecture of
B. and M. Shapiro.
This is equivalent to the following result:

If all critical points of a
rational function lie on a circle
in the Riemann sphere (for example, on the real line),
then the function maps this circle
into a circle.

\vglue12pt
\section{Introduction}

Two rational functions $f_1$ and $f_2$ will be called equivalent if
$f_1=\ell\circ f_2$, where $\ell$ is a fractional-linear transformation.
Equivalent rational functions have the same critical sets.
\specialnumber{1}\proclaim{Theorem}
\label{thm1}
If all critical points of a rational function $f$ are real{\rm ,} then~$f$ is
equivalent to a real rational function{\rm .}
\endproclaim

Lisa Goldberg \cite{Goldberg} addressed the following question:
how many equivalence classes of rational functions of degree $d$
with a given critical set of $2d-2$ points may exist?
She reduced this to the following
problem of enumerative geometry:
\nonumproclaim{Problem P} Given $2d-2$ lines in general
position in projective space ${\Bbb C\Bbb P}^d${\rm ,} how many projective
subspaces of codimension $2$ intersect all of them\/{\rm ?}
\endproclaim

To explain this reduction, to each rational function
$$f(z)=\frac{a_0+\ldots+a_dz^d}{b_0+\ldots+b_dz^d}$$
of degree $d$,
we associate a projective subspace $H_f\subset \CP^d$
defined by
the following system of two equations in homogeneous coordinates
\begin{eqnarray*}
&a_0x_0+\ldots+a_dx_d=0,\\
&b_0x_0+\ldots+b_dx_d=0.
\end{eqnarray*}
This gives a bijective correspondence between the set of equivalence classes
of rational functions of degree $d$ and the set of
subspaces of $\CP^d$ of codimension $2$ which do not intersect the
image of the {\it rational normal curve}
$E(z)=\break (1:z:\ldots:z^d)\in\CP^d,\; z\in\CP^1$.
By straightforward computation one can verify 
that $z_0$ is a critical point of $f$ if and only if the tangent line at
$E(z_0)$
to the rational normal curve intersects the
subspace $H_f$.

The answer to Problem P, going back to Schubert \cite{Schubert}
(see also \cite{Kleiman, Kleiman2}), is
\begin{equation}
\label{catalan}
u_d=\frac{1}{d}\binom{2d-2}{d-1},\quad\mbox{the $d^{\rm th}$ Catalan number.}
\end{equation}
Thus the result is

\specialnumber{A} \proclaimtitle{Goldberg \cite{Goldberg}}
\proclaim{Theorem} The number of equivalence
classes
of rational functions of degree $d$ with given $2d-2$ critical points
is at most $u_d${\rm .}
\endproclaim
\vglue-16pt
 We prove
\vglue-16pt
\phantom{yoursocool}
\specialnumber{2}\proclaim{Theorem}\label{thm2} For given $2d-2$ distinct real points{\rm ,} there exist at least
$u_d$ classes of
real rational functions of degree $d$ with these critical
points{\rm .}
\endproclaim

Theorems A and 2 imply Theorem 1. 

In general, even if the lines in Problem~P are real,
the subspaces of codimension $2$ might not be real \cite{Kleiman}.
Fulton asked the following
general question (see \cite[p.\ 55]{Fulton}):
how many solutions of real equations can be real, particularly for enumerative
problems?
We refer to a recent survey \cite{Sottilesurv} of results related
to this question.
A specific conjecture for
the Problem P was made by Boris and Michael Shapiro
(see, for example, \cite{Sottile3}):
if the lines in question are tangent
to the rational normal curve at real points,
 then all $u_d$ solutions
of the problem are real. Our Theorem 2 implies that this conjecture is true.

To reformulate Theorem 1, we write a rational function  as a ratio of
polynomials without a common factor, $f=f_1/f_0$, and suppose for simplicity that
$\infty$ is not a critical point of $f$. Then critical points of $f$ coincide
with zeros of the Wronski determinant
$W(f_0,f_1)=f_0f_1^\prime-f_0^\prime f_1,$ and Theorem 1 is equivalent to
the following: {\it 
if the Wronskian of two polynomials has only real zeros{\rm ,}
then these polynomials can be made real by a linear transformation with constant
coefficients.} A more general conjecture of B. and M. Shapiro states that this
is true for any number of polynomials. It is not enough to require that
the Wronskian has real coefficients. Indeed, if $f_1(z)=z^3+3iz^2$ and
$f_0(z)=z-i$, then $W(f_0,f_1)=2z^3+6z$ has real coefficients, but no nontrivial
linear combination of $f_0$ and $f_1$ is a real polynomial.

A general discussion of the B. and M. Shapiro conjectures, with
experimental evidence and bibliography,
is contained in \cite{Sottile2}, \cite{Sottile3}.
For the related problem of pole
assignment
in the theory of automatic control
we refer to\break \cite{BB}, \cite{B}.

As a corollary from his main result in \cite{Sottile}, Sottile proved
that there exists an open (in the usual
topology)
set $X\subset\R^{2d-2}$, such that for $x\in X$ there exist $u_d$ classes of
real rational
functions of degree $d$, whose critical set is given by $x$.
Theorem 2 was proved by Sottile for $d=3$, and tested,
using computers, for $d\leq 9$.
The computation for $d=9$ ($u_9=1,430$) is due to
Verschelde \cite{Versch}.

It is interesting that
our proof of Theorem 1 is based on the fact that two different
enumerative problems have the same sequence of integers as their solution.
These two problems are Problem P and the one in Lemma 1 below.
We prove Theorem 2 in Sections 2--6 and derive
Theorem 1 in Section 7.

The scheme of our proof of Theorem 2 is following. We consider the unit circle $\T$
instead of the real line. Let $R$ be the set of rational
functions of degree $d$, mapping $\T$ into itself,
having $2d-2$ distinct critical points
in $\T$, and being
properly normalized.
For $f\in R$ we introduce a ``net'' $\gamma(f)=f^{-1}(\T)$, considered modulo
symmetric (with respect to $\T$)
normalized homeomorphisms of the Riemann sphere, preserving orientation.
A net partitions the Riemann sphere into simply-connected regions; each
of these regions is mapped by $f$ homeomorphically onto a component of $\bC\backslash\T$.
Equivalence classes of nets are combinatorial objects,
describing topological properties of rational functions $f\in R$.
To describe a function $f\in R$ completely, we need one more piece of data,
which we call a labeling. It is a function on the set of edges of a net, which
assigns to each edge the length of its image.
We give a precise description of all nets $\gamma$ (modulo equivalence)
and labelings which may occur. It is important that, for a fixed $\gamma$, the space
of possible labelings has simple topological structure: it is a convex polytope.
To recover $f$ from its labeled net, we first construct a ramified covering
$g:\bC\to\bC$, which maps each edge of $\gamma$ homeomorphically onto an arc
of the unit circle, whose length is specified by the label of this edge.
Furthermore $g$ maps each component of the complement $\bC\backslash\gamma$
homeomorphically onto an appropriate component of $\bC\backslash\T$.
Once such ramified covering $g$ is constructed, the Uniformization theorem
implies the existence of a homeomorphism $\phi:\bC\to\bC$, such that $f=g\circ\phi^{-1}$.

This construction leads to a parametrization of the set $R$
by equivalence classes of labeled nets.
Similar parametrizations for
polynomials and trigonometric polynomials were studied by Arnold in
\cite{Arnold1}, \cite{Arnold2}, and for meromorphic functions on arbitrary Riemann surfaces
by Vinberg \cite{Vinberg}, who used the nets.
The dual graph of a net of a meromorphic function
is known in classical function theory
as a ``line complex,'' or a Speiser graph \cite{GO}, \cite{Wittich}. It is
essentially our tree $S$, which will be described in Section 2.

Nonequivalent nets correspond to nonequivalent rational functions.
For a fixed net $\gamma$,
each labeling defines a rational function of the class $R$.
Taking the critical set of this rational function, we obtain  a map
$\Phi$ from the space of labelings of $\gamma$
to the space of critical sets on the unit circle. We prove that $\Phi$ is
surjective.
So for a given critical set, each $\gamma$ gives a
rational function of our class $R$, and it remains to count all possible classes of nets
$\gamma$. It turns out that there are exactly $u_d$ of them (Lemma 1).

The main difficulty is the proof of surjectivity of $\Phi$.
This is achieved by a version of the ``continuity method'' going back to
Poincar\'{e} and Koebe (see, for example, \cite[Ch.\ V, \S6]{Goluzin}),
but we have to use different
tools from topology. We extend $\Phi$ to a map between closed polytopes
and show that the extended map is continuous (Sections 3 and 4). This is done
using a normal families argument, Lemma \ref{lemma2}.
An analysis of the boundary behavior of $\Phi$ in Section 5
permits us to prove surjectivity
using a topological argument in Section 6.

We thank M. Bonk, F. Nazarov and B. Shapiro 
for stimulating discussions.
Bonk suggested the subtle normalization N5 and (\ref{nor2}), which
makes our argument in Section 4 work. We also thank the referee for
valuable remarks.

We prove Theorem 1 only for $d\geq 3$, because
it is trivial for $d=2$, and because our proof would require a modification in
this case.

We fix an integer $d\geq 3$.
{}The map $\s:\bC\to\bC,\; \s(z)=1/{\bar z}$ will be called
the {\it symmetry}.
A map will be called {\it symmetric} if it commutes with the symmetry.
A set will be called {\it symmetric} if the symmetry leaves it invariant.
All homeomorphisms and ramified coverings
of the Riemann sphere $\bC=\CP^1$, except the involutions like $\s$,
are assumed to preserve orientation.
For a region $D$ we denote by $\partial D$ its oriented
boundary (so that the region is on the left). The unit circle $\T$
is always oriented anticlockwise, so $\T=\partial \U$, where $\U$ is the unit disc.
The words ``distance'' and ``diameter'' refer to
the spherical Riemannian metric on the Riemann sphere. It is obtained from the
standard embedding of $\bC$ as the unit sphere in $\R^3$.

\section{Nets, labelings and critical sequences}
\advance\eqcount by 1
 
A {\it cellular decomposition} of a set $X\subset\bC$
is a finite partition of $X$ into sets, called cells, each of them
homeomorphic to an open unit disc $\U^k\subset \R^k,\,
k=0,1,2;$ (by definition, $\U^0=\{\mbox{one point}\}$), and has closure
homeomorphic to the closed disc $\overline{\U}^k$.
The cells are called vertices,
edges and faces, according to their dimension.
The {\it degree} of a vertex is the number of edges to whose boundaries this vertex belongs.
A {\it net} $\gamma\subset\bC$
is the union of edges and vertices of some cellular decomposition of $\bC$,
which satisfies conditions N1--N5 below.
\begin{itemize}
\item[N1.] $\gamma$ is symmetric, that is $\s(\gamma)=\gamma$.
\vspace{.1in}

\item[N2.] $\T\subset\gamma$.

\item[N3.] There are $2d-2$ vertices, all belong to $\T$ and have degree $4$.
 
\item[N4.] The point $1\in\T$ is a vertex.
\end{itemize}

A cellular decomposition which satisfies N1--N4 is completely determined by its
net $\gamma$, so we permit ourselves to speak of vertices, edges and faces
{\it of a net}. Because of N3, each face $G$ has an even number of boundary vertices.
For every $\gamma$ satisfying N1--N4 we choose certain distinguished elements
as follows.
Let $v_0=1,$ and $v_1$ be the next vertex anticlockwise on $\T$.
There is a unique face $G_0$ in the unit disc,
whose boundary contains at least four vertices, $v_0$ and $v_1$ among them.
Let $v_{-1}$ be the vertex preceding $v_0$ on $\partial G_0$.
So when tracing $\partial G_0$ according to its orientation, we consecutively
encounter $v_{-1},v_0,v_1$ in this order. We also introduce two edges on the boundary
of $G_0$: $e_1=[v_0,v_1]$ and $e_{-1}=[v_{-1},v_0]$. One of these two edges,
$e'$ belongs to $\T$, the other, $e''$ does not.
Thus we have double notation for these two edges. For every $\gamma$ satisfying N1--N4
there is a unique choice of the {\it distinguished elements}
$G_0,v_{-1},v_0,v_1,e_{-1},e_1,e',$ and
$e''$ (see Figure 1). The vertices of $\gamma$ will be enumerated as $v_1,\ldots,v_{2d-2}$,
anticlockwise on $\T$, so that
$v_{2d-2}=v_0$, and $v_{-1}=v_N,$ for some $N=N(\gamma)\in [3,2d-3]$.
Our last assumption about nets is the normalization condition
\begin{itemize}
\item[N5.] $v_{-1}=e^{-2\pi i/3},\quad v_0=1,\quad\mbox{and}\quad v_1
=e^{2\pi i/3},\quad\mbox{the cubic roots of $1$}.$
\end{itemize}
(The particular choice of these three points on $\T$ is irrelevant).
Two nets $\gamma_1$ and $\gamma_2$ are called {\it equivalent}
if there exists a symmetric homeomorphism $h$ of the sphere $\bC$, such
that $h(\gamma_1)=\gamma_2$,
and $h$ leaves each cubic root of $1$ fixed.
Such $h$ induces a bijective correspondence
between the cells of the corresponding cellular decompositions,
so we can speak of a vertex, an edge or a face of a class of nets.
Each distinguished element described
above is mapped by $h$ onto a distinguished element with the same name.
We denote by $[\gamma]$ the equivalence class of a net $\gamma$.

For a net $\gamma$ we denote by $V,E$ and $Q$ the sets of its
vertices, edges and faces, respectively.
Euler's formula implies $|Q|=2d$ and $|E|=4d-4$.
We denote by $Q_\U\subset Q$ the subset of faces which belong to $\U$,
and by $E_\T$ the subset of edges, which belong to $\T$.

Figure 1 shows all nets for $d=4$ with distinguished faces and vertices.
For aesthetic reasons we ignored N5 in this picture.
\figin{gabri1}{450}
\centerline{Figure 1. All nets for $d=4$. Only the parts in $\bU$ are shown.}

\specialnumber{1}\proclaim{Lemma}\label{prop2} There exist exactly
$u_d$
classes of nets with $2d-2$ vertices{\rm ,} where $u_d$ is the Catalan number
{\rm (\ref{catalan})}{\rm .} 
\endproclaim

This can be found in \cite[Exercise 6.19 {\bf n}]{Stanley}.
This exercise contains 66 combinatorial problems with the Catalan
numbers as the answer (see also Exersise 6.25 for algebraic interpretations
of these numbers).
Stanley uses notation $C_n=u_{n+1}$.

To each net $\gamma$ corresponds the {\it dual graph}
$S$ of the cellular
decomposition
of $\bU$ defined by $\gamma$. More precisely, each
vertex $q=q_G$ of $S$ corresponds to a face $G=G_q\in Q_\U$, and two vertices
of $S$ are connected by an edge $\tau=\tau_e$ in $S$ if the two corresponding
faces in $Q_\U$ have a common edge $e=e_\tau$ in $\gamma\cap\U$.

Let $\hat S$ be the graph obtained by the following extension of $S$:
for every edge $e\in E_\T$, a vertex $q_e$ and an edge $\tau_e$
connecting $q_e$ with $q_G$ are added to $S$, where $G$ is the
face in $Q_\U$ with $e\in\partial G$.

It is easy to see that $S$ and $\hat S$ are trees. We designate $q_0=q_{G_0}$
to be the root of these trees. Notice that the edges of $\hat S$ are
in bijective correspondence with edges of $\gamma$ in $\bU$,
and the edges of $S$ correspond to the edges of $\gamma$ in $\U$.
There is a natural partial order on the vertices of a rooted tree, so that
the root is the unique minimal element. Thus the tree $S$ defines a
partial order on faces in $Q_\U$:
\begin{equation}
\label{order}
G'<G\quad\mbox{if the path in $S$ from $q_0$ to $q_G$ passes through $q_{G'}$}.
\end{equation}
We can also
order the set of faces in $Q_\U$ into a sequence
$G_0,\ldots,G_{d-1}$ so that, for every $k\in [1,d-1]$,
the face $G_k$ has exactly
one common boundary edge with the union of the faces $G_0,\ldots,G_{k-1}$.
Such ordering is always compatible with the partial order (\ref{order}):
\begin{equation}
\label{mo}
\mbox{for every}\quad m,n\in [0,d-1]\quad G_n<G_m\quad\mbox{implies}\quad n<m.
\end{equation}
We will use repeatedly the possibility of such ordering.

For a net we define a function $\sigma:Q\to\{1,
-1\}$, called
{\it parity}. We put $\sigma(G_0)=1,$ for the distinguished face,
and then $\sigma(G)\sigma(G')=-1$
if the faces $G$ and $G'$ have a common edge on their boundaries.
Such parity function exists for every cellular decomposition whose
vertices have even degree. With our normalization
$\sigma(G_0)=1$, the parity function is unique.
 
A {\it labeling} of a net is a nonnegative function on the set of edges,
$p:E\to\R$, satisfying the following conditions:
\begin{equation}
\label{symme}
p(\s(e))=p(e)\quad\mbox{for every}\quad e\in E,
\end{equation}
where $\s$ is the symmetry with respect to $\T$,
\begin{equation}\label{cond}
\sum_{ e\subset\partial G}p(e)=2\pi\quad\mbox{for every}\quad G\in Q_\U,
\end{equation}
and
\begin{equation}
\label{normlab}
p(e_1)=p(e_{-1})=2\pi/3.
\end{equation}
A pair $(\gamma,p)$ is called a {\it labeled net}.
Two labeled nets $(\gamma_1,p_1)$ and $(\gamma_2,p_2)$ are
{\it equivalent} if there exists a symmetric homeomorphism
$h:\bC\to\bC$, fixing the three cubic roots of $1$, and having the properties
$h(\gamma_1)=\gamma_2$, and
$p_2(h(e))=p_1(e)$ for every edge $e$ of $\gamma_1$.
 
A labeling $p$ is called {\it degenerate} if $p(e)=0$ for some edges $e\in E$,
otherwise it is called {\it nondegenerate}. The space of all
labelings is a closed convex polytope $\overline{L}_\gamma$ in
the affine subspace $A$ of $\R^{4d-4}$
defined by (\ref{symme}), (\ref{cond}) and (\ref{normlab}).
Its interior $L_\gamma$ with respect to $A$, the set of nondegenerate
labelings, is homeomorphic to a cell of dimension $2d-5$.

The statement about dimension will not be used, but it can be verified
in the following way.
First, using the equations (\ref{symme}), we
eliminate the variables $p(e)$ for all edges in $\bC\backslash\bU$.
The number of remaining variables is $3d-3$.
Each of the equations (\ref{cond}) corresponds to
a face $G\in Q_\U$. This face $G$ has at least one boundary edge on $\T$,
which does not belong to the boundaries of other faces in $Q_\U$.
Thus each equation in (\ref{cond}) contains a variable which does not
show in other equations. So the codimension of the affine subspace
defined by all equations (\ref{cond}) is $d$. Equations (\ref{normlab})
increase the codimension to $d+2$. So the dimension of $L_\gamma$ is
$2d-5$.

A {\it critical sequence} corresponding to $\gamma$
is a map $c:V\to\T$, which leaves $v_0,\,v_1$ and $v_N=v_{-1}$ fixed,
and preserves the (nonstrict) cyclic order.
We describe critical sequences by
nonnegative functions $l:E_\T\to\R$,
For $k=1,\ldots,2d-2$, the value $l([v_{k-1},v_k])$
is defined as
the length of the arc
$[c(v_{k-1}),c(v_{k})]$,
of $\T$, described anticlockwise from $c(v_{k-1})$ to $c(v_{k})$.
This function $l$ has the following properties:
\begin{equation}
\label{l1}
l([v_0,v_1])=\sum_{k=2}^{N}l([v_{k-1},v_k])=
\sum_{k=N+1}^{2d-2}l([v_{k-1},v_k])=\frac{2\pi}{3},\quad
\mbox{and}\quad l\geq 0,\hskip.15in
\end{equation}
where $N=N(\gamma)$. Thus we identify the set of all critical sequences
with the convex polytope $\overline{\Sigma}_\gamma$, described by
(\ref{l1}).
This polytope is a product of two simplexes of dimensions $N-2$ and $2d-N-3$,
so its dimension is $2d-5$. The interior $\Sigma_\gamma$ of our
polytope consists of
critical sequences with $l>0$. We call such critical sequences
{\it nondegenerate}, and the critical sequences in $\overline{\Sigma}_\gamma
\backslash\Sigma_\gamma$ {\it degenerate}.

We denote by $R^*$ the class of all
rational functions of degree at most $d$, which preserve the unit circle,
whose critical points all belong to the
unit circle, and which satisfy the normalization condition
\begin{equation}\label{nor2}
f(z)=z,\quad f'(z)=0,\quad\mbox{for}\quad z\in\{1,e^{2\pi i/3},e^{-2\pi i/3}\}.
\end{equation}
This normalization implies that two different functions of the class $R^*$ are
never equivalent.

For each class of nets $[\gamma]$, we consider a subclass
$R_\gamma\subset R^*$ defined by the following
condition:
\begin{equation}\label{nor1}
f^{-1}(\T)\in[\gamma].
\end{equation}
It follows from (\ref{nor1}) that $R_\gamma$ consists of rational
functions of degree $d$ with simple critical points,
which coincide with the vertices of $f^{-1}(\T)$.
Furthermore, (\ref{nor2}) and (\ref{nor1}) imply that $f$
maps the distinguished face $G_0$ of the net $f^{-1}(\T)$ onto the unit disc.

It will follow from the results of Section 3 that $R_\gamma\neq\emptyset$
for every $\gamma$.

\section{Construction of a map}
\advance\eqcount by 9
\begin{equation}\label{map}
F_\gamma:\overline{L}_{\gamma}\rightarrow R^*\times\overline{\Sigma}_{\gamma}.
\end{equation}
In this section, for each net $\gamma$, we construct a map (\ref{map}), where
$\overline{L}_{\gamma}, R^*$ and $\overline{\Sigma}_{\gamma}$
were introduced in Section 2,
with the following
properties:
\begin{equation}
\label{prop1}
F_\gamma(L_\gamma)\subset R_\gamma\times\Sigma_\gamma.
\end{equation}
If $p$ is a nondegenerate labeling, and $(f,c)=F_\gamma(p)$, \pagebreak 
then $c$ is the sequence of critical points of $f$.
An additional property, related to the boundary behavior of $F_\gamma$,
is stated in Proposition \ref{147} below.
In Section 4 we will prove that the second component $\Phi_\gamma$
of $F_\gamma$
is continuous, and in
Section 6 that $\Phi_\gamma$
is surjective.

To construct our map $F_\gamma$, we fix $\gamma$ and
a labeling $p\in\overline{L}_\gamma.$
We introduce the following notation.
Let $Z$ be the union of edges $e$ with $p(e)=0$,
and $D$ the component
of $\bC\backslash Z$, containing $G_0$.
We claim that
\begin{equation}
\label{noloops}
0<p(e)<2\pi\quad\mbox{for every}\quad e\subset D.
\end{equation}
The left inequality follows immediately from the definition of $D$.
To prove the right inequality, we suppose without loss of generality
that $e\subset\bU$, and use the tree $S$ introduced
in Section 2. Let $G\subset D$ be a face in $Q_\U$ whose boundary
contains $e$.
Then there is a path in $S$ from the root $q_0$ to $q_G$.
It is easy to see that all $G_q$ for $q$ in this path belong to
$D$.
The labels of all edges along this path are positive, because the whole
path belongs to $D$. It follows from (\ref{cond}) that the labels
of all edges of this path are less than $2\pi$. Thus no edge in $\partial G$
can have label $2\pi$.

We put
$B=\bC\backslash D$ and introduce an equivalence relation in $\bC:\;
x\sim y$ if $x$ and $y$ belong to the same component of $B$.
Let $Y=\bC/\sim$ be the factor space, and $w:\bC\to Y$ the projection map.

Since $D$ is connected, every component of $B$ is contractible, hence $Y$
is a topological sphere, so we can identify it with the Riemann
sphere. The symmetry $\s:\bC\to\bC$ is an involution which
leaves every point of $\T$ fixed. Since every component of $B$ contains
a vertex, it intersects $\T$.
It follows that each component of $B$ is symmetric.
So $Y$ also has an involution, such that $w$ splits the involutions. This means that
the identification of $Y$ with $\bC$ can be made in such a way that
\begin{equation}\label{wi}
w:\bC\to Y\cong\bC,\quad w(x)=w(y)\quad\mbox{if and only if}\quad x\sim y 
\end{equation}
is symmetric. In particular $w(\T)=\T$. Furthermore, in view of (\ref{normlab}),
no component of $B$
can contain two cubic roots of unity, so
we can arrange that
$w(v)=v$, for each cubic root $v$ of $1$.
The cellular decomposition of $\bC$ defined
by $\gamma$ generates via $w$ a cellular decomposition $X=X(p)$ of $Y$, so
that the cells of $X$ are $w(C),$ where $C$ are the cells of the
original decomposition. If the labeling $p$ is nondegenerate,
then $w$ is a homeomorphism.

We are going to construct a continuous map
$g^*:\overline{D} \to\SS\cong\bC$, where $\SS$ is another copy
of the Riemann sphere.
As a first step of our construction of $g^*$,
we define a continuous map $\tilde g:\gamma\cap\overline{D}\to\T\subset\SS$.
To do this, we orient
the edges
of $\gamma$ in the following way. Each edge $e\in E$ belongs to the boundaries
of exactly two faces; let $G$ be that one with $\sigma(G)=1$.
Then $e\subset\partial G$ by definition inherits positive orientation of $\partial G$.

We are going to define $\tilde g$, so that the following condition be
satisfied for every $e\subset\overline{D}$:
\begin{eqnarray}\label{C}
\nonumber
&\mbox{\it if $p(e)>0$, then $\tilde g$
maps $e$ onto an arc of $\T$ of length $p(e)$,}\\
&\mbox{\it homeomorphically, respecting orientation,}\\
\nonumber
&\mbox{\it and if $p(e)=0$, then $\tilde g$ maps $e$ into a point.}
\end{eqnarray}
In particular the edges in $\partial D$ are mapped into points, but closures of
all edges in $D$ are mapped homeomorphically onto their images.
This follows from (\ref{noloops}).

First we define $\tilde g$ on $\partial G_0$, so that condition (\ref{C}) is satisfied,
and $\tilde g(v_0)=1$.
Condition (\ref{cond}) with $G=G_0$ ensures that there is a unique way
to define such continuous $\tilde g$ on
$\partial G_0$. Furthermore, (\ref{normlab}) implies that $\tilde g$ fixes
all three cubic roots of $1$.

Now we order all faces of $\gamma$ in $D\cap \U$ into a sequence
$(G_0,G_1,\ldots,G_m)$ so that for every $k=1,\ldots,m$ the face
$G_k$ has exactly one common boundary edge $e^*$ with
\begin{equation}\label{collection}
\bigcup_{j=0}^{k-1}\partial G_j.
\end{equation}
The existence of such ordering was explained in Section 2, before (\ref{mo}).

Suppose that $\tilde g$ has been already defined on
the edges in (\ref{collection}). In particular, it is defined on
the edge  $e^*\in\partial G_k$.
Condition (\ref{cond}) with $G=G_k$ allows
us to extend $\tilde g$ to all other edges
in $\partial G_k$, so that (\ref{C})
is satisfied.

After $\tilde g$ is defined for all edges in $\bU$,
we extend it to the edges in\break
$\overline{D}\backslash\bU$
by symmetry.
This construction defines a symmetric continuous map\break
$\tilde g:\gamma\cap\overline{D}\to~\T$,
which sends every component of $\partial D$ to a point.

Notice that for every face $G$, the map $\tilde g:\partial G\to\T$
has degree $\pm 1$, and is monotone; that is $\tilde g|_{\partial G}$
preserves or reverses the nonstrict cyclic order.
As a next step, for each face $G\subset D$, we extend $\tilde g$
to a continuous map
$g^*:\overline{G}\to\overline{\U}\subset\SS$, if $\sigma(G)=1$,
or $g^*:\overline{G}\to\SS\backslash \U$, if $\sigma(G)=-1$,
so that the restriction on $G$ is a homeomorphism onto the image. This can
be done for every continuous monotone map $\partial G\to\T$ of degree $\pm 1$.

It is clear, that this extension of $\tilde g$ into the interior of
components $G\in Q$, $G\subset D,$ can be made symmetrically;
that is
\begin{equation}\label{sym}
g^*\circ\s=\s\circ g^*.
\end{equation}
Finally we extend $g^*$ to a continuous map $\bC\to\SS$ so that it is constant
on every component of the set $B=\bC\backslash D$. Then $g^*(x)=g^*(y)$ whenever $x\sim y$, the equivalence
relation $\sim$ in (\ref{wi}). It follows that
$g^*$ factors as $g^*=g\circ w,$ where $w$ is the continuous map in (\ref{wi}).
Here $g$ is a continuous map $Y\to\SS$.
\pagebreak

If $C$ is a cell of the cellular decomposition defined by $\gamma$, then $w$ and $g^*$
map $C$ in the same way: either homeomorphically or to a point. It follows that
$g$ maps every closed
cell of the form $w(\overline{C})$ homeomorphically onto the image. Furthermore, the
cells $w(C)$ make a cellular
decomposition $X$ of $Y$, so $g$ is a ramified
covering. It can be ramified only at the vertices of $X$.
If the labeling $p$ is nondegenerate, that is $w$ in (\ref{wi}) is a
homeomorphism, all vertices of $X$ have order $4$, and $g$ is ramified
exactly at these vertices, having local degree $2$ at each vertex.

There exists a unique conformal structure on $Y$, which makes $g$
holomorphic.
By the Uniformization theorem \cite{Ahlfors}, \cite{Goluzin}, there exists a unique homeomorphism
$\phi:Y\to\bC$, normalized by
\begin{equation}\label{no}
\phi(e^{-2\pi i/3})=e^{-2\pi i/3},\quad \phi(1)=1,\quad
\phi(e^{2\pi i/3})=e^{2\pi i/3},
\end{equation}
and such that $f=g\circ\phi^{-1}$ is a holomorphic map $\bC\to\SS$, that is
a rational function. It is easy to see that $f$ is nonconstant and has
degree at most $d$.

This function is the first component of $F_\gamma(p)$ in (\ref{map}).
The second component is
\begin{equation}
\label{777}
c:v\mapsto\phi\circ w(v),\quad v\in V,
\end{equation}
which is a critical sequence in $\overline{\Sigma}_{\gamma}$. Indeed,
by the symmetry property (\ref{sym}) and the symmetry of
the normalization (\ref{no}), $\phi$ is symmetric. Applying
(\ref{sym}) again, we conclude that
our rational function $f$ is symmetric, and that all values of the function
$c$ belong to $\T$.
An important consequence of our construction of $F_\gamma$ is the following
proposition, which we state in terms of function $l$ as in (\ref{l1}):
\specialnumber{1}\proclaim{Proposition}\label{147}
Let $\Phi_\gamma:\overline{L}_\gamma\to\overline{\Sigma}_\gamma$
be the second component of
the map $F_\gamma${\rm ,} and $l=\Phi_\gamma(p)$ for some $p\in \overline{L}_\gamma${\rm .}
Then $l(e)\neq 0$ if and only if $e\subset D\cap\T${\rm .}
\endproclaim

This follows from (\ref{777}), taking into account that $\phi$ is a
homeomorphism, and $w$ collapses exactly those edges of $\gamma$ which do not
belong to $D$.

If the labeling $p$ is nondegenerate, then the map $w$ in (\ref{wi}),
(\ref{777}) is a
homeomorphism, which implies that all $c(v),\, v\in V$, are distinct
and coincide
with critical points of $f$. In this case we have $l(e)>0$ for
all edges $e\in E_\T$. So the second
component $\Phi_\gamma$ of $F_\gamma$
maps the set of nondegenerate labelings $L_\gamma$ to the set
of nondegenerate critical sequences $\Sigma_\gamma$, and $\gamma$ is equivalent
to $f^{-1}(\T)$ via $\phi\circ w$,
and we have (\ref{prop1}).

Now we show that our map $F_\gamma$ in (\ref{map}) is well defined,
that is $f$ and $c$ are independent of the choice of a labeled net
within its equivalence class, and also independent of the extensions of $g$
into the interiors of the components $G$.
This independence follows from \pagebreak

\specialnumber{2}\proclaim{Lemma}\label{lemma1}
Let $X_i,\; i=0,1,$ be cell complexes{\rm ,} $h'$ a bijection between their cells such that
$h'(\partial C)=\partial h'(C)${\rm ,}
$Y$ a topological space and $f_i:X_i\to Y$ two continuous maps{\rm ,} whose restrictions
to every closed cell are homeomorphisms onto the image{\rm ,} and
$f_1(C)=f_0(h'(C))$ for every cell $C$ in $X_1${\rm .} Then there exists a homeomorphism
$h$ such that $f_1=f_0\circ h${\rm .}
\endproclaim 

\demo{Proof} We define $h$ on every cell $C$ in $X_1$
as $f_{0,h'(C)}^{-1}\circ f_1|_{C}$, where $f_{0,h'(C)}^{-1}$
is the inverse of the restriction $f_0|_{h'(C)}:h'(C)\to f_0(h'(C))$.
\enddemo

Applying Lemma \ref{lemma1} to two rational functions $f_0$ and $f_1$, constructed
from equivalent labeled nets, we conclude that $f_0=f_1\circ h$, where $h$ is a
homeomorphism of the Riemann sphere. This homeomorphism is evidently conformal
and fixes
three points; thus $h=\id$ and $f_1=f_0$.

\specialnumber{3}\proclaim{Lemma}\label{wd}
The critical sequence $c=\Phi_\gamma(p)$ is well defined{\rm ,} that is it
depends only on the class of labeled nets $([\gamma],p)${\rm .}
\endproclaim

\demo{Proof} Consider the cellular decomposition $X$, introduced after equation
(\ref{wi}). If $v$ is a vertex of $X$
of degree at least $4$, then $z=\phi(v)$
is a critical point of $f$, so $c(v)$ is well defined.
Suppose now that $v^1,\ldots,v^m$ is a maximal
chain of vertices of $X$ of degree $2$, which means that there are edges in $X$ between
these vertices, but no other edges connecting $v^1$ or $v^m$ to vertices of degree $2$.
There is a unique way to extend this chain by adding
$v^0$ and $v^{m+1}$, vertices of degree at least $4$, so that $v^0$ is connected
to $v_1$ and $v^m$ to $v^{m+1}$ by edges of $X$. Then $z^j=\phi(v^j),\;
j\in\{ 0,m+1\},$ are
critical points of $f$, and $a_j=f(z^j),\; j\in\{ 0,m+1\},$ corresponding critical values.
The restriction of $f$ onto the arc $[z^0,z^{m+1}]\subset\T$ maps this arc
homeomorphically onto the arc $[a_0,a_{m+1}]\subset\T$. Then the position of the points
$z^j=\phi(v_j),\; j=1,\ldots,m,$ is determined from the fact that the length of each
arc $[f(z^k),f(z^{k+1})]\subset\T$ is equal to $p(w^{-1}([v^k,v^{k+1}]))$, the label of
an edge of $\gamma$.
\enddemo

\vglue-12pt
\section{Continuity of $\Phi$}
 \advance\eqcount by 18

For a fixed $\gamma$,
the second component of our map $F_{\gamma}$ in (\ref{map}) is a map between
two closed polytopes
\begin{equation}\label{phis}
\Phi:\overline{L}\to\overline{\Sigma},
\end{equation}
where $L=L_{\gamma},\,\Sigma=\Sigma_{\gamma}$ and $\Phi=\Phi_\gamma$.
In this section we
prove that $\Phi$ is continuous.

Suppose that $p_1\in\overline{L}$; we are going to prove that $\Phi$ is continuous at $p_1$.
Let $p_0$ be a point close to $p_1$.
Using the notation, similar to that
introduced in Section 3, before (\ref{noloops}),
we consider the sets $Z_i$ and the regions $D_i,\; i=0,1.$
In addition, let $B_1,\ldots,B_m$ be the complete
list of components of $\partial D_1$.
Then $B_j\subset Z_1$ for $j=1,\ldots,m$.
If $p_0$ is close enough to $p_1$, we may assume
\begin{equation}\label{as}
Z_0\subset Z_1,\quad\mbox{and thus}\quad D_1\subset D_0.
\end{equation}
We have the maps $w_i:\bC\to Y_i\cong\bC,\; i=0,1$ as in (\ref{wi}), and
$f_i=g_i\circ\phi^{-1}_i,$
$i=0,1$, defined in Section 4.
All maps involved in our argument are shown on the diagram below, where we use
double notation $\Bbb S=\bC$ as in Section 3.
For every vertex $v\in V$ the ``critical value''
$a(v)=g_1\circ w_1(v)$ is defined,
\begin{equation}
\label{!!}
a_k=a(v_k)=f_1\circ\phi_1\circ w_1(v_k)\in{\Bbb S},\; k=1,\ldots,2d-2.
\end{equation}
The actual set of critical values of $f_1$ is a subset of
$\{ a_k\}$ which might be proper.
\figin{xy1}{1000}
 
We choose arbitrary $\delta>0$. Then there exists $\varepsilon>0$, such that
the open discs $U_k$ of radii $\varepsilon$
around $a_k$ are either disjoint or coincide,
and have the property that every component $K$
of the preimage of their union under $f_1$ has diameter less than $\delta$.
We set $\varepsilon_1=\varepsilon/(8d)$ and
suppose that
\begin{equation}\label{close}
|p_1(e)-p_0(e)|<\varepsilon_1\quad\mbox{for every}\quad e\in E.
\end{equation}
In particular, in view of (\ref{as}),
\begin{equation}
\label{small}
|p_0(e)|<\varepsilon_1\quad\mbox{for}\quad e\subset \bigcup_{j=1}^mB_j.
\end{equation}
The set
\begin{equation}\label{set}
H={\Bbb S}\backslash\bigcup_{k=1}^{2d-2}U_k
\end{equation}
has a cell decomposition with two $2$-dimensional cells $C$ and $C^*$, where
$$\overline{C}=\bU\backslash\bigcup_{k=1}^{2d-2}U_k,$$ and $C^*$ the symmetric
cell to $C$.
We choose $1$-dimensional cells of this decomposition to be arcs of the unit circle and
arcs of the circles $\partial U_k$, and for $0$-dimensional cells
we take the points of intersections of the circles $\partial U_k$
with the unit circle.

Let $H_1$ be the preimage of  the set $H$ in (\ref{set}) under $f_1$.
Then $H_1$ has a cell
decomposition, which is the preimage of our cell decomposition of (\ref{set}),
and $f_1|_{H_1}$  maps every cell of this decomposition homeomorphically.
In fact $f_1|_{H_1}$ is a covering, because $f_1$ has no
critical points
in $H_1$.

It follows from (\ref{small}) that
\begin{equation}\label{diam}
\diam f_0\circ\phi_0\circ w_0(B_j)<\varepsilon/2,\quad\mbox{for each component
$B_j$ of $\partial D_1$},\hskip.5in
\end{equation}
because $B_j$ is made of at most $4d-4$ edges $e$, whose labels $p_0(e)$
are at most $\varepsilon_1$
each.
Furthermore, (\ref{close}) and (\ref{!!}) imply that
\begin{equation}\label{position}
\dist\left( f_0\circ\phi_0\circ w_0(v_k),a_k\right)<\varepsilon/2,
\quad k=1,\ldots,2d-2,
\end{equation}
so $f_0$ has no critical values in $H$.
As every $B_j$ contains at least one vertex $v_k$ of $\gamma$, we conclude
from (\ref{diam}) and (\ref{position}) that
\begin{equation}
\label{onemore}
f_0\circ\phi_0\circ w_0(B_j)\cap H=\emptyset,\quad j=1,\ldots,m,
\end{equation}
where $H$ was defined in (\ref{set}).
Let $H_0$ be the component of $f^{-1}_0(H)$ intersecting
$\phi_0\circ w_0(G_0)$.
As $f_0$ has no critical values in $H$,
the restriction $f_0|_{H_0}:H_0\to H$ is a covering. It follows from
(\ref{onemore})
that $H_0\subset\phi_0\circ w_0(D_1)$, because every component
$B_k$ of $\partial D_1$
is mapped by $f_0\circ\phi_0\circ w_0$ into $\bC\backslash H$.

The cell decomposition
of $H$ defined above pulls back to $H_0$, and $f_0$ maps each closed cell of this
pullback onto a cell in $H$ homeomorphically.
Notice that each open cell of $H_0$ is contained in a unique cell of the form
$\phi_0\circ w_0(C)$ for some cell $C\subset D_1$ of $\gamma$. A similar
statement
holds for cells of $H_1$. This defines a bijection between cells of $H_1$ and those of
$H_0$ which commutes with the boundary operator $\partial$.
So Lemma \ref{lemma1} can be applied to $f_i|_{H_i}$. We
conclude that
\begin{equation}\label{main}
f_1=f_0\circ h\quad\mbox{on}\quad H_1,
\end{equation}
where $h:H_1\to H_0$ is a homeomorphism. Evidently $h$ is holomorphic, and its boundary
values on $\partial H_1$ belong to $\partial H_0$. Moreover, the components
of $\partial H_1$ separating $H_1$ from the cubic roots of $1$, are
mapped to components of $\partial H_0$ separating the same cubic roots of $1$
from $H_0$.
Now we use the following

\specialnumber{4}\proclaim{Lemma}\label{lemma2}
Suppose that a finite set $X=\{ x_1,x_2,\ldots,x_n\}\subset\bC${\rm ,} is given{\rm ,} such that
$n\geq 3,$ and $x_1,x_2$ and $x_3$ are the cubic roots of $1${\rm .}
Then for every $\eta>0$ there exists $\delta\in(0,\eta)$ with the following property{\rm .}
Let $J_1,\ldots,J_n$ be disjoint open Jordan regions of diameter less than $\delta${\rm ,}
$x_k\in J_k,\; k=1,\ldots,n${\rm ,} and $h$ be an injective holomorphic function
$$h:\bC\backslash\bigcup_{k=1}^nJ_k\to\bC,$$
such that for $k\leq 3$ the curves $h(\partial J_k)$ separate $x_k$ from the two other
cubic roots of $1${\rm .} Then
$$\dist(h(z),z)<\eta,\quad\mbox{whenever}\quad \dist(z,X)\geq\eta.$$
\endproclaim

\vglue-9pt
{\it Proof} (compare \cite[Theorem 13]{Belinski}). Our proof is by contradiction. Suppose
that there is a sequence $(\delta_j)\to 0$ and a sequence $(h_j)$, which satisfies all
conditions, but
\begin{equation}\label{contra}
\dist(h_j(z_j),z_j)\geq \eta
\end{equation}
for some $\eta>0$ and some points $z_j$ with $\dist(z_j,X)\geq\eta.$ It is easy to see
that the closed domains $R_j=\bC\backslash\cup_{k=1}^nJ_{k,j}$ tend to
$\bC\backslash\{ x_1,\ldots,x_n\}$ and that all functions $h_j$ omit three
cubic roots of $1$ in their domains. By Montel's criterion \cite{Montel},
\cite{Goluzin},
$(h_j)$ is a normal family and we can select
a convergent subsequence. The limit $h$ of this subsequence is a holomorphic injective
function $h$ in $\bC\backslash\{ x_1,\ldots,x_n\}$, which omits the three cubic roots
of $1$. By the Great Theorem of Picard all points $x_k$ are removable singularities,
so $h$ extends to a fractional-linear map. But this fractional-linear map also fixes
three points, the cubic roots of $1$, so it is the identity. This contradicts (\ref{contra}).
\hfill\qed\vglue12pt

Applying Lemma \ref{lemma2} to $h:H_1\to\bC$, we obtain that
$\dist(h(z),z)<\eta,\; z\in H_1$, so
the critical sequences $(\phi_1\circ w_1(v)),\; v\in V$ and $(\phi_0\circ w_0(v)),\; v\in V$
are $\eta$-close.
So our map (\ref{phis}) is continuous.

 \vglue-9pt
\section{Boundary behavior of $\Phi$}
\advance\eqcount by 29

Our goal is to prove that
$\Phi_\gamma:L_\gamma\to\Sigma_\gamma$
is surjective. This will be achieved with the help of Lemma \ref{polytopes}
in Section 6. To verify that the conditions of this lemma are satisfied,
we need to show that the preimages of closed faces of
$\overline{\Sigma}_\gamma$ are homologically trivial.
These preimages can be complicated, so we begin
with an analysis of preimages of open faces.
In this section, a net $\gamma$ is fixed, so we do not show explicitly
the dependence of various objects on $\gamma$.

Suppose that a convex polytope $K$ is described by
$${\mathbf A}x=b,\quad x\geq 0,\quad x\in\R^n,$$
where ${\mathbf A}$ is an $m\times n$ matrix, and $b\in\R^m$.
An {\it open face} of $K$ is defined as
\begin{equation}
\label{qva}
A_W=\{ x\in K:x_j>0,\;j\in W\;\mbox{and}\;x_j=0,\;
j\notin W\},
\end{equation}
where $W$ is a subset of $\{1,\ldots,n\}$, such that
$A_W\neq\emptyset$.
A {\it closed face} is the closure of an open face, and
\begin{equation}
\label{closface}
\overline{A}_W=\{ x\in K:x_j=0,\;
j\notin W\}.
\end{equation}
Vertices are open faces and closed faces simultaneously.
All open and closed faces are nonempty convex sets.
We also notice that each closed face is a finite union of open faces.

We are going to
apply these definitions to the convex polytope $\overline{\Sigma}$, described by
(\ref{l1}), in the space of real-valued functions $l$ on $E_\T$.
First we state precisely which subsets $W\subset E_\T$ define
open faces, that is for which $W$ the set (\ref{qva}) is nonempty.
It follows from (\ref{l1}) that the necessary and
sufficient conditions are: $[v_0,v_1]\in W$, and that each of the two sequences
$$[v_1,v_2],\ldots,[v_{N-1},v_N]\quad\mbox{and}\quad
[v_N,v_{N+1}],\ldots,[v_{2d-3},v_{2d-2}]
$$
contains at least one edge in $W$.
Using the notation $e',e''$ for the distinguished edges,
introduced in Section 2, these conditions can be restated as
\vglue4pt
 (a)  $e'\in W,$ and
\vglue4pt(b)  there is at least one edge in $W\backslash\{ e'\}$ on each side of $e''$.
\vglue4pt

For a set $W\subset E_\T$, satisfying these conditions, we define the
open face $A_W$ of $\overline{\Sigma}$ by
\begin{equation}
\label{*}
A_W=\{ l:l(e)>0\;\mbox{for}\; e\in W\;\mbox{and}\; l(e)=0\;\mbox{for}\;
e\in E_\T\backslash W\}.
\end{equation}

We introduce a partial order on the set of
open faces: $A_1\prec A_2$ if and only if $\overline{A}_1\subset\overline{A}_2.$
From our definition (\ref{*}) it follows  that
\begin{equation}
\label{monotone}
{A}_{W_1}\prec{A}_{W_2}\quad
\mbox{if and only if}\quad W_1\subset W_2.
\end{equation}

To characterize the preimage $\Phi^{-1}(A_W)$, of an open face, we use
Proposition 1 from Section 3 and Lemma \ref{face2} below. To state
this lemma, we need the following
notation, similar to that used in Sections 3 and 4.
For $p\in\overline{L}$, we define
$Z(p)$ as the union of the closed edges $e$ in $E$
such that $p(e)=0$, and $D(p)$ as the connected component
of $\bC\backslash Z(p)$ containing
$G_0$. Notice that $D(p)$ always contains at least three boundary edges of $G_0$,
including $e_{-1}$ and~$e_1$.
This follows from (\ref{cond}) with
$G=G_0$ and (\ref{normlab}).
Let
\begin{equation}
\label{E^0}
E^0(p)=\{ e\in E:e\subset\partial D(p)\cap\bU\}
\end{equation}
and
\begin{equation}
\label{E(p)}
E(p)=\{ e\in E:e\subset D(p)\cap\T\}.
\end{equation}

Figure 2 shows the part in $\bU$
of a net $\gamma$ with $d=5$, the set $E^0(p)$ (bold lines), and the set $E(p)$,
which consists of the edges
$[v_0,v_1],[v_2,v_3],[v_5,v_6]$,
$[v_6,v_7]$ and
$[v_7,v_8]$ on $\T$.

\figin{gabri2}{500}
\centerline{Figure 2: $S(p)$ in dotted lines, $E^0(p)$ in bold.}
\vglue9pt\noindent 
It is clear from the definitions (\ref{E^0}) and (\ref{E(p)}), that
the set $E^0(p)$ determines $E(p)$ uniquely.
The opposite is also true: 
\specialnumber{5}\proclaim{Lemma}\label{face2}
For $p\in\overline{L}${\rm ,} the set $E(p)\subset E${\rm ,}
uniquely determines the set $E^0(p)${\rm .}
\endproclaim

\demo{Proof}
We use the rooted tree $\hat S$, introduced in Section 2.
Let $S(p)$ be the subtree of $\hat S$ spanned by $q_0$ and
$\{ q_e:e\in E(p)\}$.
Figure 2 shows
the tree $S(p)$ in dotted lines.

We claim that $S(p)=S'(p)$, where $S'(p)$ is the subtree of $\hat S$ spanned by
$\{ q_x:x\subset D(p)\cap\bU\}$ (here $x$ may stand for a face or an edge of
$\gamma$).

By definition of $E(p)$, we have $S(p)\subset S'(p)$.
It remains to prove $S'(p)\subset S(p)$,
which means that
$D(p)\cap\bU$ contains exactly those faces $G\in Q_\U$
which have the property $q_G\in S(p)$.
Thus $E(p)$ uniquely determines $D(p)$,
and $D(p)$ uniquely determines $E^0(p)$.

To prove our claim, suppose that $S'(p)\not\subset S(p)$.
Since both $S(p)$ and $S'(p)$ belong to the tree $\hat S$,
there exists a leaf $q$ of $S'(p)$ which does not belong
to $S(p)$.
If $q=q_e$, where $e\in E_\T$, then $e\subset D(p)$,
hence $q_e$ is a vertex of $S(p)$, in contradiction to our choice of $q$.
Suppose now that $q=q_G$, where $G\subset D(p)$, is a face in $Q_\U$.
Let $S_q$ be the path in $S$ connecting $q$ and $q_0$.
Conditions (\ref{cond}) imply that
$0<p(e)<2\pi$ for every $e=e_\tau,\;
\tau\in S_q$ (See also (\ref{noloops}).)
This implies that there is an edge $e\subset\partial G$ such that
$\tau_{e}\notin S_q$ and $0<p({e})<2\pi$.
Since $q\notin S(p)$, we have $\tau_e\notin S(p)$.
If $e\subset\T$, we have a contradiction with the definition of $S(p)$.
Otherwise, the other face in $Q_\U$,
having the edge $e$ on its boundary,
belongs to $D(p)$, and $G$ is not a leaf of $S'(p)$,
again a contradiction.

\hfill$\Box$
\vspace{.1in}

It follows from Proposition 1 of Section 3 and
(\ref{E(p)}) that
the preimage of an open face $A_W\subset\overline{\Sigma}$ is
\begin{equation}
\label{pqr}
\Phi^{-1}(A_W)=\{ p\in\overline{L}:
E(p)=W\}.
\end{equation}
By Lemma \ref{face2}, for every $p\in\overline{L}$, the set $W=E(p)$
uniquely determines a set $E^0_W=E^0(p)$.
According to the remark before Lemma \ref{face2}, $E^0(p)$
uniquely determines $E(p)$.
Thus (\ref{pqr}) can be rewritten as
\begin{equation}
\label{preip}
\Phi^{-1}(A_W)=\{ p\in\overline{L}:p(e)=0\;\mbox{for all}\; e\in E^0_W,
\;\mbox{and}\; p(e)>0\;\mbox{for
all}\; e\in W\}.
\end{equation}

Now we prove that these preimages (\ref{preip})
are nonempty.

\specialnumber{6}\proclaim{Lemma}\label{snakes2}
For each subset $W\subset E_\T$ satisfying
{\rm (a)} and {\rm (b)} in the beginning of this section{\rm ,}
there exists $p\in\overline{L}_\gamma$ such that $W=E(p)${\rm .}
\endproclaim

\demo{Proof}
Given a subset $W$ of edges of $\gamma\cap\T$ satisfying
{\rm (a)} and {\rm (b)}, let us define a subtree $S_W$
of the tree $\hat S$,
as the union of all paths in $\hat S$ connecting vertices $q_e$, for $e\in W$,
with $q_0$.
The labeling $p$ is defined inductively along the tree $\hat S$,
starting from the vertex $q_0$.
As $W$ contains
at least one edge, other than $e'$, at each side of $e''$,
we have $\tau_{e''}\in S_W$, and there is at least
one edge $e$ of $G_0$, other than $e'$ and $e''$,
such that $\tau_e\in S_W$.
Let $m\ge 1$ be the number of all such edges.
We define $p(e)=2\pi/(3 m)$ for each of them, $p(e')=p(e'')=2\pi/3$,
and $p(e)=0$ for all other edges of $G_0$.
This guarantees that (\ref{cond}) is satisfied\break for $G_0$.
Notice that $0<p(e)<2\pi$ for an edge $e\subset\partial G_0$
if and only if $\tau_e\in S_W$.

Suppose now that the values of $p(e)$ are defined for
all edges of faces $G_q\in Q_\U$, with
$q$ in a subtree $S'$
of $S$ containing $q_0$, so that $0<p(e)<2\pi$
if and only if $\tau_e$ belongs to $S_W$, and (\ref{cond}) is satisfied.
If $S'=S$, then the labeling $p$ is complete.
Otherwise, there exists a vertex $q^*$
in $S\setminus S'$ which is an extremity
of an edge $\tau^*$ of $S$ with another
extremity of $\tau^*$ being in $S'$.
Let $G^*=G_{q^*}$ and $e^*=e_{\tau^*}$.
Since an extremity of $\tau^*$ belongs to $S'$,
the label $p(e^*)$ is already defined.

If $p(e^*)=2\pi$ or $p(e^*)=0$, then $\tau^*$ does not belong
to $S_W$; hence all other boundary edges of $G^*$ do not belong to $S_W$.
In the first case, we define $p(e)=0$ for all edges $e\ne e^*$ of $G^*$.
In the second case, we choose an edge $e^{**}\ne e^*$ of $G^*$
and define $p(e^{**})=2\pi$ and $p(e)=0$ for all other edges of $G^*$.
Then (\ref{cond}) is satisfied for $G=G^*$.

If $0<p(e^*)<2\pi$, then $\tau^*$ belongs to $S_W$.
Since $e^*\notin\T$, there is at least one other
edge $e$ of $G^*$ such that $\tau_e$ belongs to $S_W$.
Let $n\ge 1$ be the number of all such edges.
We define $p(e)=(2\pi-p(e^*))/n$ for all these edges,
and $p(e)=0$ for all other edges $e\ne e^*$ of $G^*$.
Again we have (\ref{cond}) for $G=G^*$.
\pagebreak

Now the values of $p(e)$ are defined for
all edges of faces $G_q\in Q_\U$,
for the vertices $q$ of a connected subtree $S''$
of $S$ obtained by adding $\tau^*$ and $q^*$ to $S'$,
which concludes our inductive step. We extend our labeling $p$
to edges in $\bC\backslash\U$ by symmetry, so that (\ref{symme})
is satisfied.
The labeling $p$ constructed in
this way satisfies (\ref{cond}), (\ref{normlab}) and $W=E(p)$.
\enddemo

The closure of the set (\ref{preip}) is
\begin{equation}
\label{closure}
\overline{\Phi^{-1}(A_W)}
=\{ p\in\overline{L}:p(e)=0\;\mbox{for all}\; e\in E^0_W\},\end{equation}
which is nonempty and convex. Actually $\overline{\Phi^{-1}(A_W)}$ is
a closed face of $\overline{L}$.

Now we begin a study the intersection pattern of these sets (\ref{closure}),
which will be continued in Lemma \ref{newlemma2}.

\specialnumber{7}\proclaim{Lemma}\label{newlemma1} If $A_1\succ A_2\succ\ldots\succ A_k$
is a decreasing chain of open faces of $\overline{\Sigma}${\rm ,}
then the intersection of closures
of their
preimages
$\overline{\Phi^{-1}(A_1)}\cap\ldots\cap\overline{\Phi^{-1}(A_k)}$ is a
nonempty convex subset of $\overline{L}$ {\rm (}\/a closed face\/{\rm ).}
\endproclaim

\demo{Proof} The sets $\overline{\Phi^{-1}(A_j)}$ are convex,
so their intersection is convex. It remains to verify
that the intersection is nonempty.

We have $A_j=A_{W_j}$ for some $W_j\subset E_\T$.
As in the proof of Lemma \ref{snakes2}, for each $j\in[1,k]$,
we define a subtree $S_j\subset\hat S$ 
as the union of all paths in $\hat S$ connecting vertices $q_e,\, e\in W_j$,
with
the root $q_0$. We also define $E^0_j$
as the set of all edges $e$ of $\gamma$, such that $\tau_e\in\hat S\backslash
S_j$, and $\tau_e$ has a vertex in $S_j$.
It is easy to check that these definitions are consistent with notation
of Lemmas \ref{face2} and \ref{snakes2}: if $W_j=E(p)$ then $E^0_j=E^0(p)$.

We have the following inclusions:
$$W_1\supset\ldots\supset W_k,\quad\mbox{and}
\quad S_1\supset\ldots\supset S_k.$$
The first inclusion follows from the assumption of the Lemma and
(\ref{monotone}), the second follows from the first one, and the definition
of $S_j$. We assume without
loss of generality that $A_1=\Sigma$, so $W_1=E_\T$, and $S_1=\hat S$.
According to (\ref{closure}),
$$\overline{\Phi^{-1}(A_j)}=\{ p:p(e)=0\;\mbox{for all}\; e\in E^0_j\}.$$
Thus we have to show that there exists a labeling $p$, such that
$$p(e)=0\quad\mbox{for all}\quad e\in\bigcup_{j=1}^kE^0_j.$$

To construct this labeling $p$, we order
the set of faces in $Q_\U$ into a sequence
$G_0,\ldots,G_{d-1}$, such that for every $n\in [1,d-1]$,
the face $G_n$ has exactly
one common boundary edge with the union of faces $G_0,\ldots,G_{n-1}$.
(Such ordering was explained in Section 2, before (\ref{mo}).)

First we construct $p(e)$ for the edges $e$ in $\partial G_0$, as it is
done in the proof of Lemma \ref{snakes2}, using $W_k$ as $W$.
For these edges $e$,
we have $p(e)>0$ if and only if
$\tau_e\in S_k$.
Hence $p(e)>0$ implies
$\tau_e\in S_j$, and thus $e\notin E^0_j$ for all $j\in[1,k].$

Suppose that $p$ is already defined for all boundary edges
of faces
$G_0,\ldots,$ $ G_{n-1}$, for some $n<d$, so that
\begin{equation}
\label{proty}
p(e)>0\quad\mbox{only if}\quad e\notin E^0_j,\quad\mbox{for all}\quad
j\in[1,k].
\end{equation}
We want to extend $p$ to the boundary edges of $G=G_n$, so that
the property (\ref{proty}) is preserved.

Let $m\in[1,k]$ be the integer, such that $q_G\in S_m\backslash
S_{m+1}$ ($S_{k+1}:=\emptyset$).
Consider the path $\Gamma$ in the tree $S$ from $q_0$ to $q_G$. Let $q_{G'}$
be the vertex on this path preceding $q_G$. Then $G'<G$ in the sense
of the partial order defined in (\ref{order}), and this implies by (\ref{mo})
that $G'\in\{ G_0,\ldots,G_{n-1}\}$.
There exists exactly one boundary edge $e^*$ in $\partial G\cap\partial G'$.
This is the only edge in $\partial G$, on which $p$ is defined so far.
Since $q_G\in S_m$, we have $q_{G'}\in\Gamma\subset S_m$.
This implies $\tau_{e^*}\in S_m$.
There is at least one more boundary edge $e^{**}$ of $G$, such that
$\tau_{e^{**}}\in S_m$. This is because all leaves of $S_m$ are in $\T$;
hence $q_G$ is not a leaf.
We define $p(e^{**})=2\pi-p(e^*),$
and $p(e)=0$ for all edges on $\partial G$, other than $e^*$ and $e^{**}$.

Notice that on this inductive step, the only new edge for which a positive
value of $p$ was defined is the edge $e^{**}$. Now we are going to prove
that the condition (\ref{proty}) was preserved on the inductive step.
Since $q_{e^{**}}\in S_m$, we have $q_{e^{**}}\in S_j$ and $e^{**}\notin
E^0_j$ for $j\leq m$.
Since $q_G\notin S_{m+1}$, $e^{**}$ does not belong to any $E^0_j$ for
$j>m$. Otherwise, the face $G''\neq G$ such that $e^{**}\subset\partial G
\cap\partial G''$ would have the property $q_{G''}\in S_j$;
hence the path from $q_{G''}$ to
$q_0$, which contains $q_G$, would belong to $S_j$,
which is impossible since $q_G\notin S_j$.

This inductive procedure defines $p(e)$ for all edges. The labeling
$p$ we defined satisfies (\ref{proty}), as required.
\enddemo

\vglue-12pt
\section{Surjectivity of $\Phi$ and proof of Theorem 2}
\advance\eqcount by 39
\vglue-9pt
 
In this section we use homology groups with integral coefficients.
We call a topological space {\it homologically trivial} if it has
the same homology groups as one point. In particular, such set is nonempty
and connected.
Nonempty convex sets are homologically trivial. Thus
Lemma \ref{newlemma1} of the previous section shows that our map $\Phi$
satisfies the conditions of the following

\specialnumber{8}\proclaim{Lemma}
\label{newlemma2}
Let $\Phi:\overline{L}\to\overline{\Sigma}$ be a continuous map of
closed polytopes{\rm ,} such that for every $k\geq 1$ and for every
 decreasing chain
$A_1\succ\ldots\succ A_k$ of open faces of $\overline{\Sigma}${\rm ,} the set
$$\overline{\Phi^{-1}(A_1)}\cap\ldots\cap\overline{\Phi^{-1}(A_k)}$$
is homologically trivial{\rm .} Then for every $k\geq 1$ and for every
decreasing chain
$A_1\succ\ldots\succ A_k$ of open faces of $\overline{\Sigma}${\rm ,} the set

$$\overline{\Phi^{-1}(A_1)}\cap\ldots\cap\overline{\Phi^{-1}(A_{k-1})}
\cap\Phi^{-1}(\overline{A_k})$$
is homologically trivial{\rm .} In particular{\rm ,} $\Phi^{-1}(\overline{A})$ is
homologically trivial for every open face $A${\rm .}
\endproclaim

To prove this lemma we need the following result from
\cite{Godement}, Corollaire de Th\'eor\`eme 5.2.4 (of Leray).
Suppose that a compact topological space $K$ has
a finite covering by its closed subsets $\{ K_j\}$, such that
all intersections
\begin{equation}
\label{pp}
K_{j_1}\cap\ldots\cap K_{j_m}
\end{equation}
 are either empty
or homologically trivial. The {\it nerve} of such a covering is defined
as the simplicial complex, whose vertices are $K_j$ and a subset of
vertices $\{ K_{j_1}\ldots,K_{j_m}\}$ defines a simplex if and only if
the intersection
(\ref{pp}) is nonempty. {\it Then $K$ has the same homology groups
as the nerve.}

Another version of this result
was proved by K. Borsuk \cite{Borsuk}. Borsuk's theorem
assumes the nonempty intersections to be absolute retracts, and
concludes that $K$ is of the same homotopy type as the nerve.
We can use either of these two results, but we prefer the homology version.
\vspace{.1in}

\demo{Proof of Lemma {\rm \ref{newlemma2}}} We use induction on $d=\dim A_k$.
For $d=0$, $A_k$ is
one point, so $\overline{A_k}=A_k$, and the assumption of the lemma
contains
its conclusion in this case.

Suppose now that a chain $A_1\succ\ldots\succ A_k$ is given,
$\dim A_k=d\geq 1$, and the conclusion of the lemma
holds for all decreasing chains whose last term is of dimension at most
$d-1$. Consider the set
$$X=\overline{\Phi^{-1}(A_1)}\cap\ldots\cap\overline{\Phi^{-1}(A_{k-1})}
\cap\Phi^{-1}(\overline{A_k})=X_0\cup\{ X_B:B\prec A_k,\;\dim B\leq d-1\},$$
where 
$$X_0=\overline{\Phi^{-1}(A_1)}\cap\ldots\cap\overline{\Phi^{-1}(A_{k})},$$
and
$$X_B=\overline{\Phi^{-1}(A_1)}\cap\ldots\cap\overline{\Phi^{-1}(A_{k-1})}
\cap\Phi^{-1}(\overline{B}),$$
for all open faces $B\prec A_k$ of dimension at most $d-1$.
Then for every collection $B_1,\ldots,B_q$ of such open faces
we have
\begin{equation}
\label{xxx}
X_0\cap X_{B_1}\cap\ldots\cap X_{B_q}=
\overline{\Phi^{-1}(A_1)}\cap\ldots\cap\overline{\Phi^{-1}(A_k)}
\cap\Phi^{-1}(\overline{B}),
\end{equation}
where $\overline{B}=\overline{B_1}\cap\ldots\cap\overline{B_q}$
is a face of dimension at most $d-1$.
The set (\ref{xxx}) is homologically trivial if and only if
$B$ is nonempty, by the assumption of induction. Similarly,
\begin{equation}
\label{xxxx}
X_{B_1}\cap\ldots\cap X_{B_q}=
\overline{\Phi^{-1}(A_1)}\cap\ldots\cap\overline{\Phi^{-1}(A_{k-1})}\cap
\Phi^{-1}(\overline{B}),
\end{equation}
where $\overline{B}=\overline{B_1}\cap\ldots\cap\overline{B_q}$
is homologically trivial if and only if $B$ is nonempty.

This means that the nerve of the closed covering $X_0\cup_B\, X_B$ of $X$
coincides with the nerve of the covering of $\overline{A_k}$ by the
closures of open faces $B\prec A_k$, including $\overline{A_k}$ itself.
By the corollary of Leray's theorem stated above,
$X$ has the same homology groups as $\overline{A_k}$,
but $\overline{A_k}$ is nonempty and convex, so $X$ is homologically trivial.
\enddemo

 From Lemma \ref{newlemma2} we conclude that the  preimages of {\it closed}
faces of $\overline{\Sigma}$ under $\Phi$ are homologically trivial,
in particular, they are nonempty and connected. 
To complete the proof of Therem 2 we use the following

\specialnumber{9}\proclaim{Lemma}\label{polytopes}
Let $\overline{L}$ be a   topological space{\rm ,}
$\overline{\Sigma}$ a finite cell complex{\rm ,} all closed cells of $\overline{\Sigma}$ are homeomorphic to
closed balls{\rm ,} and
$\Phi:\overline{L}\to\overline{\Sigma}$ a continuous mapping such that
the preimage of each closed cell in $\overline{\Sigma}$ is 
homologically trivial{\rm .}
Then $\Phi$ is surjective{\rm .}
\endproclaim

{\it Proof}.
Let $C_k(\overline{L})$ be the space of $k$-chains of $\overline{L}$
with integer coefficients.
Let $C_*$ be the corresponding chain complex,
with the natural differential
$\partial:C_k(\overline{L})\to C_{k-1}(\overline{L})$.
Let $C_*(\overline{\Sigma})$ be the corresponding
chain complex for $\overline{\Sigma}$.

For every closed $k$-cell $A$ of $\overline{\Sigma}$,
we are going to construct a chain
$W_A\in C_k(\overline{L})$ so that  $\Phi(W_A)\subset A$, $\Phi(\partial W_A)\subset \partial A$
and $\Phi^{A}_* [W_A]=[A]$.
Here $\Phi^{A}_*$ is the mapping $H_k(\Phi^{-1}(A),\Phi^{-1}(\partial A))
\to H_k(A,\partial A)$ induced by $\Phi,\;[W_A]$ is the class
of $W_A$ in $H_k(\Phi^{-1}(A),\Phi^{-1}(\partial A))$,
and $[A]$ is the class of $A$ in $ H_k(A,\partial A)\cong\Z$.

We proceed inductively on $k=\dim A$.
For $k=0$,  the preimage of the vertex $A$ is nonempty, 
and we take a point in $\Phi^{-1}(A)$ as $W_A$.
Suppose that the chains $W_A$ are defined for all cells
$A$ of $\overline{\Sigma}$ with $\dim A<k$,
so that $\Phi^{A}_* [W_A]=[A]$ and $\partial W_A = W_{\partial A}$.
Here a chain $W_C$, for a chain $C=\sum m_\nu A_\nu$, is defined as
$\sum m_\nu W_{A_\nu}$.

Let $A$ be a cell in $C_k(\overline{\Sigma})$.
Due to the induction hypothesis, $W_{\partial A}$
is a cycle in $C_{k-1}(\overline{L})$, and
$\Phi^{\partial A}_* [W_{\partial A}] = [\partial A]$.
Here $\Phi^{\partial A}_*$ is the mapping $H_{k-1}(\Phi^{-1}(\partial A))\to
H_{k-1}(\partial A)$ induced by $\Phi$.
As $\Phi^{-1}(A)$ is a homologically trivial
subcomplex of $\overline{L}$, there exists a chain
$W_A\in C_k(\overline{L})$ such that
$W_A\subset\Phi^{-1}(A)$ and $\partial W_A = W_{\partial A}$.
{}From the commutative diagram
$$\begin{array}{ccc}
H_k(\Phi^{-1}(A),\Phi^{-1}(\partial A)) &\stackrel{\partial}{\longrightarrow} &H_{k-1}(\Phi^{-1}(\partial A))\\[5pt]
 {\scriptstyle\Phi^A_*}\Big\downarrow&&   {\scriptstyle\Phi^{\partial A}_*}\Big\downarrow\\[4pt]
H_k(A,\partial A)& \stackrel{\partial}{\longrightarrow} &H_{k-1}(\partial A)
\end{array}
$$
we have $\partial \Phi^A_* [W_A] = \Phi^{\partial A}_* \partial [W_A]
= [\partial A]$.
As $H_k(A,\partial A)\cong \Z$ is generated by $[A]$ and
$\partial [A] = [\partial A]\ne 0$,
this implies $\Phi^A_* [W_A] = [A]$.

To complete the proof, we have to show that the mapping $\Phi:W_A\to A$
is surjective, for any cell $A$ of $\overline{\Sigma}$.
If this is not so, then there exists an internal point $a\in A$ not covered
by $\Phi(W_A)$.
Since $A\setminus a$ is contractible to $\partial A$, this contradicts
the condition $\Phi^A_*([W_A])=[A]\ne0$.
\hfill\qed

\demo{Proof of Theorem {\rm 2}} It remains to summarize what has been done.
In Section 3, for each net $\gamma$ we constructed a map (\ref{map}),
which transforms labelings into pairs $(f,c)$, where $f$ is a rational function
of the class $R^*$,
and $c$ a critical sequence. Restricting this map to nondegenerate
labelings we obtain rational functions of the class $R_\gamma\subset R^*$,
whose critical
points are given by the nondegenerate sequence $c$ (see (\ref{prop1})).
By Proposition 1 in Section 3, degenerate labelings produce degenerate
critical sequences, that is $\Phi^{-1}(\Sigma)=L$ in Lemma \ref{polytopes}.
This lemma implies that all
possible critical sets, consisting of $2d-2$ points,
can be obtained in this way. So each $\gamma$
produces a rational function of the class $R_\gamma$ with prescribed
critical points. We conclude by Lemma 1, that the number of rational
functions in $R^*$, with $2d-2$ prescribed critical points is at least $u_d$.
This proves Theorem 2 because, as we saw in the end of Section 2,
different functions from $R^*$ are nonequivalent.
\hfill$\Box$

\section{Proof of Theorem 1}

In  this section we derive Theorem 1 from Theorems A and 2.
The vector space $\Poly_d$ of polynomials of degree at most $d$ with
complex coefficients is identified with $\C^{d+1}$.
Every pair $(r,q)$ of nonproportional polynomials spans a $2$-dimensional
subspace in $\Poly_d$.

To parametrize the equivalence classes of rational functions of degree
$d$, we consider the Grassmannian $G(2,d+1)$,
which is the
set of all $2$-dimensional subspaces in $\Poly_d$, and the locus $D_1\subset
G(2,d+1)$ of those pairs of\break polynomials $(r,q)$, for which
$\deg r/q<d$.
Then $D_1$ is an algebraic subvariety of\break $G(2,d+1)$
of codimension $1$. Two pairs $(r_i,q_i),\; i=1,2,$ represent the same point in
$G(2,d+1)\backslash D_1$ if and only if the rational functions $r_1/q_1$ and
$r_2/q_2$ are equivalent. Thus classes of rational functions of degree $d$ are
parametrized by $G(2,d+1)\backslash D_1$.

The Wronski determinant of two nonproportional polynomials
$$W(r,q)=\left|\begin{array}{ll}r&q\\r'&q'\end{array}\right|$$
is a nonzero polynomial of degree at most $2d-2$,
whose zeros are finite critical points
of $f=r/q$, counting multiplicities, and common zeros of $r$ and $q$.
The common zeros of $r$ and $q$ are multiple zeros of
$W(r,q)$. If two pairs of polynomials define the same point in $G(2,d+1)$,
then the Wronskians of these pairs differ by a constant multiple.
The set of all nonzero polynomials of degree
at most $2d-2$, modulo proportionality,
is parametrized by $\CP^{2d-2}$.
Thus we have a regular map $\Tilde{W}:G(2,d+1)\to\CP^{2d-2},$
defined by taking the proportionality class of the Wronski determinant.

We show that $\Tilde{W}$ is a finite map \cite[p.\ 177]{AG}. This fact is known
\cite{BB}, \cite{Goldberg}, but we include a short proof.
We normalize our Wronskians, so
that the coefficient of the monomial of the smallest degree equals $1$.
Notice that each monomial
$z^n$, where $0\leq n\leq 2d-2$,
has only finitely many preimages under $\Tilde{W}$, namely
the $2$-subspaces, generated by pairs
$(z^k,z^m)$, where $k+m=n+1$ and $k\neq m$.
If $(r,q)$ represents a point in $G(2,d+1)$,
we consider the one-parametric family of points represented by
$(r_\lambda,q_\lambda),\;\lambda\in \C^*$, where
$r_\lambda(z)=r(\lambda z)$ and $q_\lambda(z)=q(\lambda z)$.
Putting $w_\lambda=\Tilde{W}(r_\lambda,q_\lambda),$
we obtain $W(r_\lambda,q_\lambda)(z)=\lambda W(r,q)(\lambda z)$,
and after normalization $w_\lambda(z)=\lambda^{-n-1}W(r,q)(\lambda z)$,
where $n\in [0,2d-2]$ is the smallest degree of monomials in $W(r,q)$.
So $\dim \Tilde{W}^{-1}(w_\lambda)=\dim\Tilde{W}^{-1}(w_1)$
for $\lambda\in\C^*$, and
$$w_0:=\lim_{\lambda\to 0}w_\lambda\quad\mbox{is}\quad w_0(z)=z^n.$$
As the dimension of preimage
is an upper semi-continuous function of the point \cite[p.\ 138]{AG},
for regular mappings
into compact spaces, that is
$$\limsup_{\lambda\to 0}\dim\Tilde{W}^{-1}(w_\lambda)\leq
\dim\Tilde{W}^{-1}(w_0)=0,$$
we conclude that $\dim\Tilde{W}^{-1}(w_1)=0$ for every $w_1\in \CP^{2d-2}$,
so the preimages are finite,
and the map $\Tilde{W}$ is finite.

Let $D_2\subset\CP^{2d-2}$ be the locus of polynomials with
multiple roots, or having smaller degree than $2d-2$.
Notice that $\Tilde{W}(D_1)\subset D_2$.
According to Theorem A, for every point $w$ in $\CP^{2d-2}\backslash D_2$
we have $|\Tilde{W}^{-1}(w)|\leq u_d$. On the other hand,
our Theorem 1 implies that for every point $w$ in
the open set ${\bf V}$ in $\RP^{2d-2}$ formed by polynomials with $2d-2$
distinct real zeros the cardinality
of $\Tilde{W}^{-1}(w)\cap G_\R(2,d+1)$ is at least $u_d$.
Here $G_\R(2,d+1)$ stands for the `real part' of the Grassmannian,
that is the collection of those $2$-dimensional subspaces which can be
generated by pairs of real polynomials. This means that
$$\Tilde{W}^{-1}({\bf V})\subset G_\R(2,d+1).$$
On the other hand, for finite maps
we have $\Tilde{W}^{-1}(w)=\lim_{w'\to w}\Tilde{W}^{-1}(w')$,
so $\Tilde{W}^{-1}(\overline{\bf V})\subset G_\R(2,d+1)$,
where $\overline{\bf V}$ is the subset of $\RP^{2d-2}$
formed by polynomials with all real zeros.
\hfill\qed

\bye